\documentclass[reqno,11pt]{amsproc}
\usepackage{amsmath, amssymb, amsthm}
\usepackage{graphicx}

\theoremstyle{plain}

\theoremstyle{definition}

\usepackage{bm}

\newcommand{\vT}{\bm{T}}
\newcommand{\ve}{\bm{e}}
\newcommand{\vth}{\bm{\theta}}

\begin{document}

\title[A numerical study of rigidity of hyperbolic splittings]{A numerical study 
of rigidity of hyperbolic splittings in simple two-dimensional maps}

\author[O.F.~Bandtlow]{O.F.~Bandtlow}
\address{
O.F.~Bandtlow\\
School of Mathematical Sciences\\
Queen Mary University of London\\
London E3 4NS\\
UK.
}
\email{o.bandtlow@qmul.ac.uk}

\author[W.~Just]{W.~Just}
\address{
W.~Just\\
Institut f\"ur Mathematik\\
Universit\"at Rostock\\
Rostock\\
Germany. 
}
\email{wolfram.just@uni-rostock.de}

\author[J.~Slipantschuk]{J.~Slipantschuk}
\address{
J.~Slipantschuk\\
Department of Mathematics\\
University of Warwick\\
Coventry\\
UK.
}
\email{julia.slipantschuk@warwick.ac.uk}

\begin{abstract}
Chaotic hyperbolic dynamical systems enjoy a 
surprising degree of rigidity, a fact which is well
known in the mathematics community but perhaps less so 
in theoretical physics circles. Low-dimensional hyperbolic systems are
either conjugate to linear automorphisms, that is,
dynamically equivalent to the Arnold cat map and its variants, or
their hyperbolic structure is not smooth.
We illustrate this dichotomy using a family of analytic maps, for which we show by means of numerical simulations that the 
corresponding hyperbolic structure is not smooth, thereby providing an example for a global mechanism which produces non-smooth phase space structures in an otherwise smooth dynamical system.
\end{abstract}

\keywords{Dynamical Systems, Chaos, Hyperbolic Maps}

\subjclass[2010]{Primary 37D20; Secondary 37E30}

\date{22 February 2024}

\maketitle

For decades, nonlinear dynamics has had a significant impact on various applications in many scientific fields,
largely due to the study of seemingly simple model systems that exhibit surprisingly complex and counterintuitive behaviour,
often with great relevance for real-world phenomena.
Evidence of this impact
can be found in the now-classic studies of one- and two-dimensional
time-discrete dynamical systems, showing universal transitions from regular 
to chaotic dynamical behaviour, via the so-called Feigenbaum and
quasiperiodic routes to chaos, which allow for fairly rigorous underpinnings through renormalisation group treatments 
\cite{Feig_JSP79, FeKaSh_PD82}.
More importantly, the relevance of the 
mathematical findings have been confirmed in a
variety of experiments (see, for example, \cite{GiMuPe_PRL81, StHeLi_PRL85}).
Furthermore, simple low-dimensional chaotic maps
are at the core of the foundations of statistical physics, as
they provide a rigorous mechanism for the emergence
of irreversibility, the approach to
thermodynamic equilibrium, and the understanding of transport properties 
in Hamiltonian time-invariant dynamical systems,
without taking recourse to phenomenological coarse graining
\cite{HaSa_PRA92, BaCo_JPA94, Dorf:99}.

While the case of one-dimensional maps has been studied thoroughly, both empirically and rigorously,
for higher-dimensional invertible maps the situation is somewhat different.
Despite a wealth of available numerical data, rigorous statements are
sparse due to the considerable technical challenges involved.
The majority of explicit calculations going beyond
numerical simulations are based on a few paradigmatic systems,
such as the Smale horseshoe or the Arnold cat map and its variants.
These studies, as well as most applied investigations, build on the concept of hyperbolicity,
see \cite[Ch.~6 and 19]{KaHa:96} for a comprehensive overview.

In order to keep technical details to a minimum we shall focus on maps
$\vT$ mapping the two-dimensional torus $[0,1)^2$ into itself, that is 
$\vT(\vth)\in [0,1)^2$ 
with $\vth=(\theta_1,\theta_2) \in [0,1)^2$.
Omitting technical details
and looking at the orientation-preserving case only, hyperbolicity
amounts to the existence of two normalised transversal vector fields
$\ve_s(\vth)$ and $\ve_u(\vth)$, the so-called
stable and unstable directions, which are invariant under the dynamics
in the sense that
\begin{equation}\label{a}
\begin{aligned}
D\vT(\vth) \ve_u(\vth)=\lambda_u(\vth) \ve_u(\vT(\vth)) \\
D\vT(\vth) \ve_s(\vth)=\lambda_s(\vth) \ve_s(\vT(\vth)) \, .
\end{aligned}
\end{equation}
Here,  $D\vT(\vth)$ denotes the Jacobian of the map, and
the so-called local contraction and expansion
rates obey
$0<\lambda_s(\vth)<1<\lambda_u(\vth)$. A diffeomorphism $\vT$ exhibiting this
property globally is referred to as an Anosov diffeomorphism.
An explicit calculation
of this hyperbolic structure, that is, closed expressions
for $\ve_u(\vth)$ and $\ve_s(\vth)$, have been obtained only in
very few cases where the Jacobian does not depend on the
phase space coordinates. These cases serve as textbook
examples for hyperbolic dynamics (see, for example, \cite{Dorf:99}).
The hyperbolic structure is a key mechanism for chaotic dynamics in
time-reversible systems, and hence a key ingredient for
irreversibility occurring in otherwise Hamiltonian dynamics. While
an explicit computation of the stable and unstable directions
poses considerable challenges, a numerical calculation is more or less 
straightforward. This can be achieved through a suitable forward or backward iteration of the linear
variational equations (\ref{a}), which mimics the simple
power method for computing leading eigenvalues and eigenvectors of
matrices. In fact, this empirical numerical scheme is also at the
heart of rigorously proving the existence of a hyperbolic splitting
via an invariant cone field using the Alekseev criterion
\cite{Alek_MUSSR68}.

We now turn to
the surprising degree of rigidity enjoyed by hyperbolic systems,
a fact well known in the mathematical dynamics community, but quite counterintuitive and
perhaps less well known in the applied sciences. Naively one would
expect that the degree of smoothness of the dynamical system is
reflected by the smoothness of the hyperbolic structure; for example, one might suspect
that analytic hyperbolic dynamical systems also have an
analytic hyperbolic splitting. This turns out to be far from the truth.
Anosov showed in \cite{A} that the hyperbolic splitting of any smooth
Anosov diffeomorphism is $C^\alpha$
(that is, $\alpha$-H\"{o}lder continuous) for some $\alpha \in (0,1)$,
but need not be $C^1$ or even Lipschitz\footnote{See \cite[p.201]{A} for an example
where the splitting is almost nowhere $C^{2/3 + \epsilon}$ for any $\epsilon > 0$.}. Here, the regularity of the hyperbolic splitting refers to that of the stable and unstable subspaces with respect to the phase space point.
Higher-order differentiability is rare, making the H\"older setting natural for
hyperbolic structures of (even smooth) Anosov diffeomorphisms.
In the special case of \emph{symplectic} smooth Anosov diffeomorphisms on the two-torus
$[0, 1)^2$, the natural setting typically becomes $C^{1+\alpha}$ with $\alpha \in (0,1)$, see, for example, \cite{Has94}.
If in this case the hyperbolic splitting is of regularity\footnote{In fact,
in this setting the hyperbolic splitting being $C^2$
implies that it is already $C^\infty$ \cite{Has92}.} $C^2$,
the diffeomorphism is in fact necessarily smoothly conjugated to a
linear Anosov diffeomorphism, that is,
it is dynamically equivalent to a variant of the cat map. See \cite{Kato_PMIHES90,Gh93,HaWi99} for these, and more nuanced related results.

This obstruction to higher regularity in arbitrarily smooth (nonlinear) hyperbolic systems is
somewhat unexpected at first sight,
and shows that
in physical dynamical systems non-smooth behaviour can be produced
by a global mechanism.
Although these results are well-established, to the best of our knowledge, it is difficult to find explicit examples in the literature where the
lack of smoothness has been demonstrated purely through numerical methods.
In this article we aim to fill this gap by explicitly demonstrating the lack of smoothness of the hyperbolic splitting for a specific class of systems that are, to a significant extent, amenable to an analytic treatment.

In order to construct the desired examples we shall resort to a class of analytic maps on the torus derived from Blaschke maps, for which full spectral information is available. This means that, in addition to ergodic invariant measures, it is also possible to compute exponential decay rates
of correlation functions for analytic observations in closed form. 
The case of one-dimensional Blaschke maps has been studied for
decades \cite{Mart_BLMS83,BePe_DS16}. However, the complete analytic treatment
of their spectral properties is a more recent development
\cite{SlBaJu_NONL13,BaJuSl_AIHP17}, for a precursor see also \cite{LeSoYu94}.
While formal generalisations to 
two dimensions have been proposed a while ago 
\cite{PuSh_ETDS08}
the construction of two-dimensional hyperbolic diffeomorphisms with
accessible spectra is a very recent finding 
\cite{SlBaJu_NONL17,PoSe_NONL23,SlBaJu22}.

\begin{figure}[!h]
\centering
\includegraphics[width=0.5\textwidth]{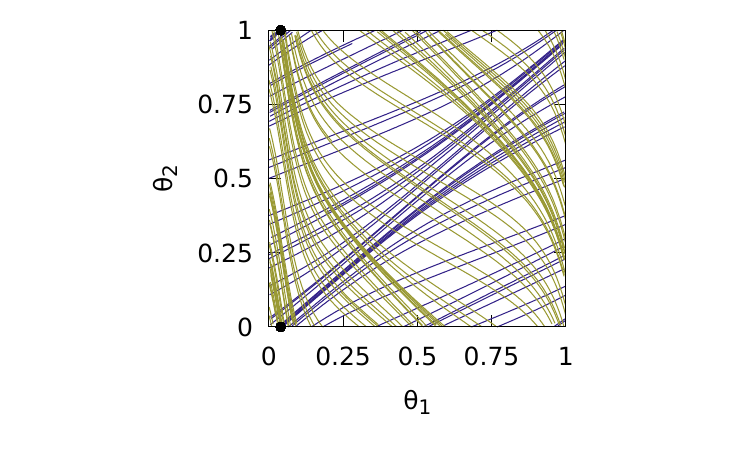}
\caption{Stable (blue) and unstable (bronze) manifold of
the fixed point (black filled circle) of the map $\vT$ given in (\ref{b}) and (\ref{c})
for $\mu=0.7$ and $\alpha=0.3$. Manifolds of finite length
have been computed numerically by forward, respectively backward, 
iteration (and an additional adaptive bisection scheme) of initial conditions 
close to the fixed point and located on the stable, respectively 
unstable, direction.
\label{figa}}
\end{figure}

More specifically, we shall investigate the properties of the hyperbolic splitting
for the following two-dimensional mixing analytic map on the torus
\begin{equation}\label{b}
\vT(\vth)=(2\theta_1+\theta_2+f(\theta_1), \theta_1+\theta_2+f(\theta_1))
\end{equation}
with
\begin{equation}\label{c}
\begin{aligned}
f(\theta)=&\frac{1}{\pi} \arctan \left( \frac{\mu \sin(2 \pi \theta-\alpha)}
{1-\mu \cos(2 \pi \theta-\alpha)} \right)  \\
 & (0\leq \mu <1, \alpha\in[-\pi,\pi)) \, .
\end{aligned}
\end{equation}
This map can be considered as a (strong) perturbation of the
Arnold cat map with (\ref{c}) denoting the deformation of the map
governed by the parameters $\mu$ and $\alpha$.
The unique fixed point of the map is given by
\begin{equation}
\vth_*=\left(\frac{1}{\pi} \arctan\left (\frac{\mu \sin\alpha}{1+\mu\cos\alpha}\right),0\right )
\end{equation}
with the Jacobian obeying $\mathop{\mathrm{det}}(D\vT(\vth_*))=1$ and 
\begin{equation} 
\mathop{\mathrm{tr}}(D\vT(\vth_*))=1+2\frac{1+\mu \cos\alpha}{1-\mu^2}\,.
\end{equation}
The above formula for the trace implies in particular that for $\mu >0$ the map $\vT$ given by (\ref{b}) is not smoothly 
conjugate to the Arnold cat map (which arises from (\ref{b}) by choosing $\mu =0$), since the Jacobians at the
unique fixed point of both maps, do not coincide.  
It turns out that $\vT$ is a symplectic map for which Lebesgue measure is invariant and mixing.
Moreover, the point spectrum of the compact Perron--Frobenius operator defined on a 
suitable anisotropic Hilbert space is given by 
$\{ (-\mu)^{|n|} \exp(i n \alpha) :  n \in \mathbb{Z}\}$, and correlation 
functions of analytic observables decay exponentially
\cite{SlBaJu_NONL17}.
Furthermore, the structure of (\ref{b}) ensures the existence
of an invariant cone field, given by the standard quadrants in the tangent space.
Hence by the Alekseev criterion \cite{Alek_MUSSR68}, $\vT$
has a hyperbolic splitting as defined in (\ref{a}).
While the hyperbolic splitting is hard, if not impossible,
to compute by analytic means,
it is fairly easy to get some impression
of the hyperbolic structure by computing the stable and unstable
manifolds of the fixed point $\vth_*$
numerically, that is, by a
forward and backward iteration of the map, see Figure~\ref{figa} (see also \cite[\S 4.3]{ChaWa}).
The hyperbolic splitting visible in this figure is not
twice continuously differentiable, since the map is not smoothly
conjugate to the Arnold cat map. However, this degree of 
roughness is difficult to discern visually.
Results in \cite{Kato_PMIHES90} provide an upper bound on the smoothness
of the hyperbolic splitting for properly nonlinear maps
and suggest that typically the splitting is, up to logarithmic
corrections, almost twice differentiable. We are going to
demonstrate this phenomenon for our map $\vT$. 
In our subsequent numerical calculations we will use the numerical
values $\mu=0.7$ and $\alpha=0.3$ for the parameters of the map.
For this setting one can on the one hand still verify hyperbolicity of the
map by an invariant cone condition, while on the other hand a non-vanishing
value for $\alpha$ avoids an additional inversion symmetry,
which would occur if $\mu \exp (i \alpha)$ were real valued. Nevertheless,
the following numerical results do not seem to depend on these
considerations\footnote{The Fortran codes used  
to produce the numerical results are available at
https://www.math.uni-rostock.de/\~{}wj/code/.}.

Since (\ref{a}) relates the local expansion and contraction rates with the stable and
unstable directions, the degree of smoothness of the hyperbolic
splitting is mirrored by the smoothness of these rates.
Hence, in order to evaluate
the degree of smoothness of the hyperbolic structure it
is sufficient to examine the smoothness of the local expansion
rate $\lambda_u(\vth)$.

\begin{figure}[!h]
\centering
\includegraphics[width=0.32\textwidth]{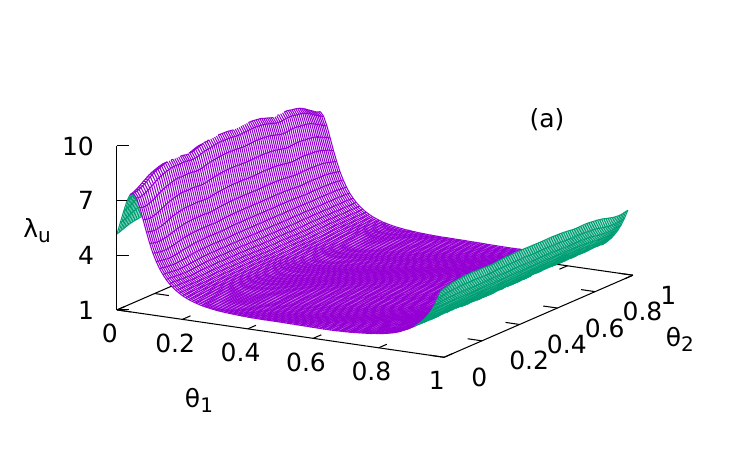}
\includegraphics[width=0.32\textwidth]{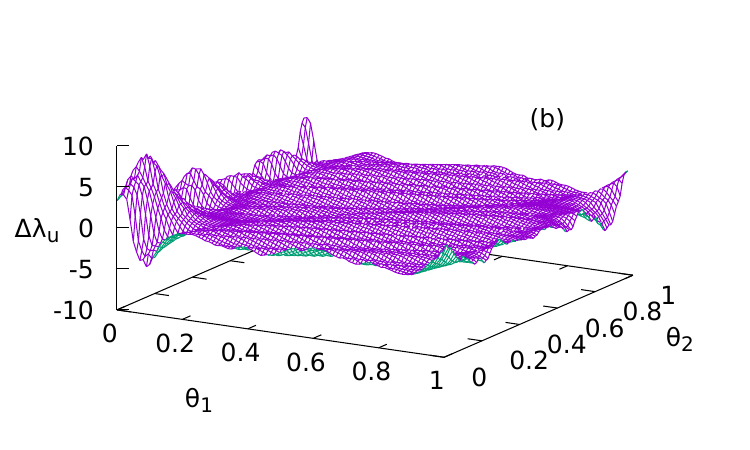}
\includegraphics[width=0.32\textwidth]{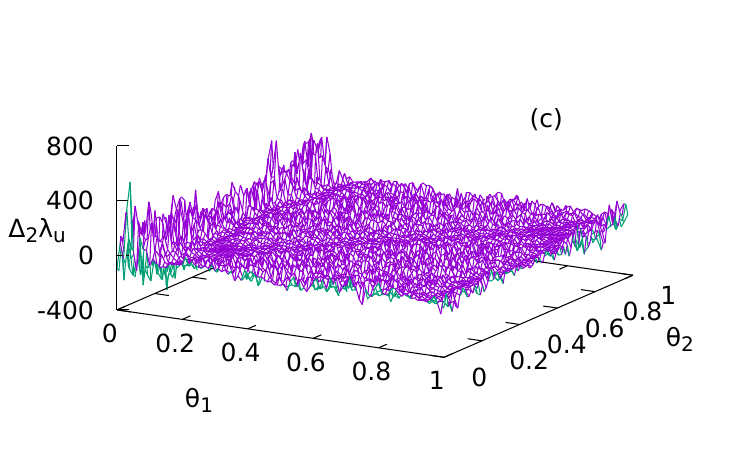}
\caption{Local expansion rates on a $100 \times 100$ square 
grid of $\vth$-values, computed from backward orbits of
length 200. (a) Local expansion rate $\lambda_u(\vth)$.
(b) First order difference quotient defined in (\ref{d}), for $h=10^{-4}$.
(c) Second order difference quotient defined in (\ref{e}), for $h=10^{-4}$.
\label{figb}}
\end{figure}

Using (\ref{a}) we numerically evaluate the local expansion rate mimicking
the power method: given a phase space point $\vth$ we compute
a backward orbit of finite length and then use forward iteration of
(\ref{a}) along this orbit to obtain a numerical approximation of
$\lambda_u(\vth)$ from the normalisation of the image vectors. 
At quadruple precision a backward orbit of length of 
about 100 is sufficient to obtain numerically accurate results 
for $\lambda_u(\vth)$ up to 30 digits. Results for the local expansion rate
are shown in Figure~\ref{figb}(a). In order to evaluate the degree
of smoothness of this graph we consider difference quotients to estimate
the derivative. Here we just show results for the partial difference quotients
in the  $\theta_2$-direction. Results for the symmetric
first order difference quotient
\begin{equation}\label{d}
\Delta \lambda_u = \frac{\lambda_u(\theta_1,\theta_2+h)-
\lambda_u(\theta_1,\theta_2-h)}{2 h}
\end{equation}
which estimates the first order partial derivative are shown in 
Figure~\ref{figb}(b), for a fixed value of the offset $h$. Again we obtain an 
apparently smooth function, indicating that the local expansion rate may
be continuously differentiable. In order to obtain evidence that the hyperbolic
splitting is not twice differentiable we shall now consider the
second order difference quotient
\begin{equation}\label{e}
\Delta_2 \lambda_u = \frac{\lambda_u(\theta_1,\theta_2+h)+
\lambda_u(\theta_1,\theta_2-h)- 2 \lambda_u(\theta_1,\theta_2) }{h^2} \, .
\end{equation}
Figure~\ref{figb}(c) shows a rough surface with large
fluctuations, suggesting that the second order difference quotient does not converge to a well-defined second 
derivative, as expected. In order to underpin this observation, namely that the hyperbolic 
splitting is not $C^2$, we shall now
consider the scaling behaviour of the difference quotients (\ref{d}) and (\ref{e}) with respect to the offset $h$ in greater detail.

\begin{figure}[!h]
\centering
\includegraphics[width=0.49\textwidth]{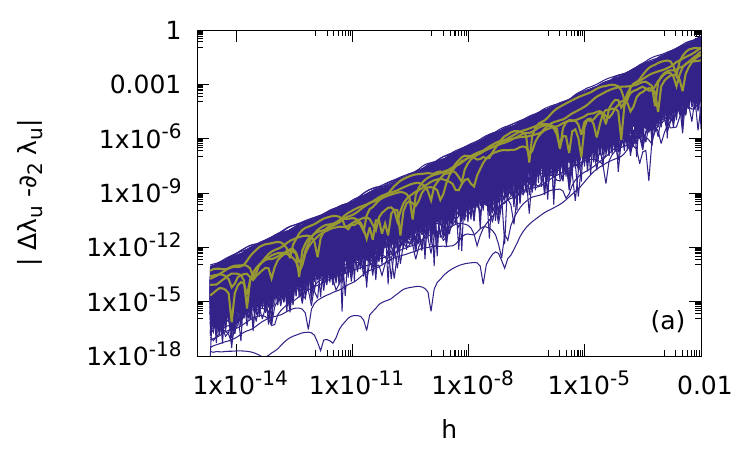}
\includegraphics[width=0.49\textwidth]{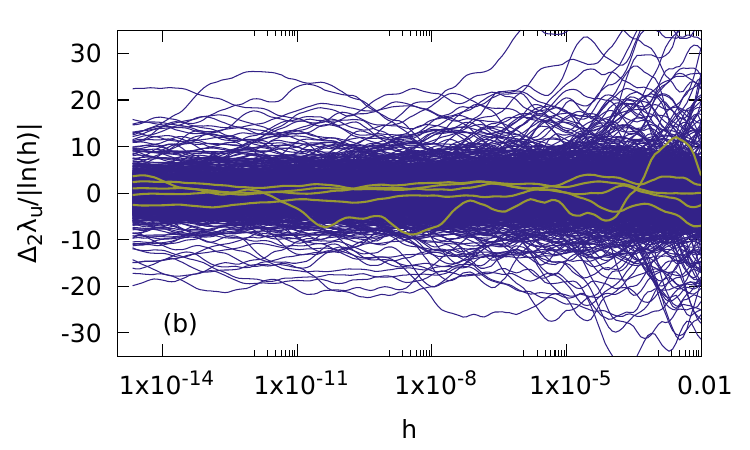}
\caption{(a) Log-log plot of the difference between $\Delta \lambda_u$ and
a numerically ``exact'' estimate for the partial derivative
(obtained via the difference quotient at $h=10^{-16}$) 
in dependence on the offset $h$. 
Data are shown for different values of $\vth$
on a  $40\times 40$ square grid (blue), and for five selected 
$\vth$-values (highlighted, bronze). 
(b) Semi-log plot of the dependence of the second order
difference quotient $\Delta_2 \lambda_u$ on the offset $h$. 
The second order difference quotient is scaled by $|\ln h|$.
Data are shown for different values of $\vth$
on a  $40\times 40$ regular grid (blue), and for five selected values
(highlighted, bronze).
\label{figc}}
\end{figure}

To begin with, we look at the convergence of the first order difference
quotient towards the limit value, that is, at the dependence of $\Delta \lambda_u$ given in (\ref{d}) 
on the offset~$h$. In order to capture possible dependencies on the phase space
coordinate we perform the numerical evaluation for $\vth$-values
on a regular grid of size $40 \times 40$. Even though one could apply more
sophisticated methods this turns out to be sufficient to estimate the
potential limit by computing the difference quotient for a very small
value of the offset. In our case, the value $h=10^{-16}$ turns out to be 
sufficient.
The dependence of the difference between $\Delta \lambda_u$ and the limit value 
is displayed in Figure~\ref{figc}(a).  Clearly one observes convergence 
of the difference quotient towards a limit at a rate of order $\mathcal{O}(h)$.
This rate is slower than the rate one would expect for an analytic 
function, but it is in line with the fact that
the local expansion rate cannot be twice differentiable. The picture
looks different for the second order difference quotient, see 
Figure~\ref{figc}(b). The second order difference quotient $\Delta_2 \lambda_u$, given in (\ref{e}), apparently
fails to converge and is even unbounded
as its behaviour is affected by logarithmic or
sub-logarithmic corrections. Even though the precise scaling is
difficult to estimate from the numerical data, the results shown in
Figure~\ref{figc}(b) suggest that 
$\Delta_2 \lambda_u=\mathcal{O}(|\ln h|)$. This matches the best known regularity result for the hyperbolic splitting in this setting, namely that it is in the Zygmund class, see \cite{Kato_PMIHES90}.
Hence the data provide compelling evidence that the hyperbolic splitting has modulus of continuity $\mathcal{O}(s | \log s |)$, implying it is
$C^{2-\varepsilon}$ for any $\varepsilon>0$, that is, it is $C^1$ and its derivative is $\alpha$-H\"older continuous for all $\alpha\in (0,1)$.

In summary, we have presented a numerical illustration of a well-established fact, namely, that, phrased in colloquial terms,
for analytic chaotic maps there occurs a dichotomy between maps with a non-smooth
hyperbolic splitting and maps that are
smoothly conjugate to the cat map or its variants. This observation points to
a mechanism creating subtle structures in dynamical systems, and
defies the naive belief that nature prefers smooth
structures. Our computations suggest
that the local expansion rate $\lambda_u(\vth)$ is continuously
differentiable but the second order derivatives fail to converge
due to logarithmic corrections.
Even though the precise form of the
logarithmic corrections is difficult to assess from numerical simulations,
the data shown in Figure~\ref{figc} are in line with rigorous results
(see, for example, \cite[Thm 3.1, Cor 3.6]{Kato_PMIHES90}),
which may not be as readily accessible to a general audience as they deserve to be. 
Furthermore, the logarithmic 
scaling
seems to occur for typical phase space points, at least in our example. Overall, the numerics suggests
that the hyperbolic splitting is $C^{2-\varepsilon}$, that is, it is
almost twice differentiable - as asserted by the theory. These fine details are certainly not
visible in simple illustrations such as Figure~\ref{figa} which superficially
indicates analytic behaviour. 

We have focussed here on the simplest case of two-dimensional symplectic
diffeomorphisms of the torus, which allowed for a straightforward
formulation of the dichotomy between conjugacy to linear diffeomorphisms
and occurrence of a non-smooth hyperbolic splitting,
and enabled a simple numerical check
in terms of a single scalar function like the local expansion rate.
It is tempting
to look at higher-dimensional cases, where, however, the dichotomy becomes
much more involved. The hyperbolic
splitting of a general Anosov diffeomorphism is a priori only H\"older
continuous, regardless of the smoothness of the map, and a numerical
check of the properties of the splitting may now require
a more laborious and numerically demanding
computation of stable and unstable eigenspaces (see, for instance,
\cite{GiChLiPo_JPA13, Noet_PD19}).
This makes the numerical detection of any kind of potential
logarithmic corrections to scaling extremely challenging,
as it requires scanning quantities over a large range of magnitudes.

We would also like to remark that the lack of higher-order regularity of
statistical quantities pertaining to a smooth dynamical system is also of
great import for applications. 
For example, the lack of differentiability of unstable Jacobians
has posed an obstacle for the efficient computation of linear response using 
Ruelle’s formula \cite{Ruel_CMP97} (see \cite{GoLi,Jia12,ChJe_NONL22,Ni23} and \cite[\S5.3]{Bal18} for background
and a
broader overview of linear response).
It is worth noting that our choice of area-preserving Blaschke maps as perturbations is not well-suited to numerically study this particular problem (due to the perturbations all sharing Lebesgue measure as their relevant invariant measure, thus 
giving rise to vanishing derivatives).

While we expect that other area-preserving perturbations of linear toral automorphisms yield comparable numerical performance,
our choice of Blaschke maps was motivated
by the fact that these allow for explicit analytic computation of
dynamical properties such as correlation decay and spectra of transfer
operators \cite{SlBaJu_NONL17, PoSe_NONL23, SlBaJu22}.
The explicit calculation of the spectra of transfer operators in these cases relies on
analytic properties of the dynamical system and the introduction of suitably
chosen anisotropic Hilbert spaces on which the transfer operator is compact 
and has a spectrum that is accessible by means of a certain composition operator.
There is no obvious relation
between these features and topological properties in phase space such as
hyperbolic splittings. In fact, we have not been able to produce a map
where a non-trivial splitting can be obtained in closed form. It is likely
that an example of this type would trigger renewed interest in
a striking global dynamical phenomenon that appears to have been settled from
a mathematical perspective.

\section*{Acknowledgements}
All authors gratefully acknowledge the support of the research presented in this article by the EPSRC grant EP/R012008/1. Additionally, J.S. acknowledges support by ERC-Advanced Grant 833802-Resonances, and
W.J. acknowledges funding by the Deutsche Forschungsgemeinschaft (DFG, 
German Research Foundation) - SFB 1270/2 - 299150580.
We would also like to thank Oliver Butterley, Fran\c{c}ois Ledrappier and Mark Pollicott for providing pointers to the mathematical literature.  

\section*{Data availability statement}

No new data were created or analysed in this study.


\vspace{1em}


\begin{thebibliography}{10}

\bibitem{Alek_MUSSR68}
V. Alekseev, {\em Quasirandom dynamical systems I: quasirandom diffeomorphisms},
Math.\ USSR Sb. {\bf 5},  73  (1968).

\bibitem{A}
D. V. Anosov,
{\em Geodesic flows on Riemann manifolds with negative curvature},
Proc.\ Steklov Inst. {\bf 90}, 818 (1967).

\bibitem{Bal18}
V. Baladi,
{\em Dynamical Zeta Functions and Dynamical Determinants for Hyperbolic Maps.
A Functional Approach.}, Ergebnisse Vol. 68, Springer, 2018.

\bibitem{BaCo_JPA94}
O.F. Bandtlow and P.V. Coveney,
{\em On the discrete time version of the Brussels formalism},
J.\ Phys. A {\bf 27},  7939  (1994).
  
\bibitem{BaJuSl_AIHP17}
O.F. Bandtlow, W. Just, and J. Slipantschuk,
{\em Spectral structure of transfer operators for expanding circle maps},
Ann.\ Inst.\ H. Poinc.\ C {\bf 34}, 31 (2017).

\bibitem{ChJe_NONL22}
N. Chandramoorthy and M. J{\'e}z{\'e}quel,
{\em Rigorous justification for the space–split
sensitivity algorithm to compute linear
response in Anosov systems},
Nonlinearity {\bf 35}, 4357 (2022). 
  
\bibitem{BePe_DS16}
C. Penrose and C. Beck, 
{\em Superstatistics of Blaschke products}, 
Dyn.\ Syst. {\bf 31} 89 (2016).

\bibitem{ChaWa}
N. Chandramoorthy and Q. Wang,
{\em An ergodic-averaging method to differentiate covariant Lyapunov vectors},
Nonlin. Dyn. {\bf 104} 4083 (2021).
  
\bibitem{Dorf:99}
J. Dorfman,
{\em An Introduction to Chaos in Nonequilibrium Statistical Mechanics}
(CUP, Cambridge, 1999).

\bibitem{Feig_JSP79}
M. Feigenbaum,
{\em The universal metric properties of nonlinear transformations},
J.\ Stat.\ Phys. {\bf 21}, 669 (1979).

\bibitem{FeKaSh_PD82}
M. Feigenbaum, L. Kadanoff, and S.J. Shenker,
{\em Quasiperiodicity in dissipative systems: a renormalization group analysis},
Physica D {\bf 5}, 370 (1982).

\bibitem{GiChLiPo_JPA13}
F. Ginelli, H. Chat{\'e}, R. Livi, and A. Politi,
{\em Covariant Lyapunov vectors},
J.\ Phys. A {\bf 46}, 254005 (2013).

\bibitem{GoLi}
S. Gou\"ezel and C. Liverani,
{\em Banach spaces adapted to Anosov systems},
Ergod.\ Th.\ \& Dyn.\ Sys. {\bf 26}, 189 (2006).

\bibitem{Gh93}
E. Ghys,
{\it Rigidit\'e diff\'erentiable des groupes fuchsiens},
Publ.\ math.\ l'IHES {\bf 78} (1993).

\bibitem{GiMuPe_PRL81}
M. Giglio, S. Musazzi, and U. Perini,
{\em Transition to chaotic behavior via a reproducible sequence of period-doubling bifurcations},
Phys.\ Rev.\ Lett. {\bf 47},  243  (1981).

\bibitem{HaSa_PRA92}
H.H. Hasegawa and W.C. Saphir,
{\em Unitarity and irreversibility in chaotic systems},
Phys.\ Rev.\ A {\bf 46},  7401  (1992).

\bibitem{Has92}
B. Hasselblatt,
{\em Bootstrapping regularity of the Anosov splitting},
Proc.\ of AMS {\bf 115} (1992). 

\bibitem{Has94}
B. Hasselblatt,
{\em Regularity of the Anosov splitting and of horospheric foliations},
 Ergod.\ Th.\ \& Dyn.\ Sys.  {\bf14}, 645 (1994). 

\bibitem{HaWi99}
B. Hasselblatt and A. Wilkinson,
{\em Prevalence of non-Lipshitz Anosov foliations},
 Ergod.\ Th.\ \& Dyn.\ Sys.  {\bf 19}, 643 (1999).

\bibitem{Jia12}
M. Jiang,
{\em Differentiating potential functions of SRB measures on hyperbolic attractors},
 Ergod.\ Th.\ \& Dyn.\ Sys. {\bf 32}, 1350 (2012). 

\bibitem{Kato_PMIHES90}
S. Hurder and A. Katok,
{\em Differentiability, rigidity and Godbillon-Vey classes for Anosov flows},
Publ.\ math.\ l'IHES {\bf 72},  5  (1990).

\bibitem{KaHa:96}
A. Katok and B. Hasselblatt,
{\em Introduction to the Modern Theory of Dynamical Systems},
(CUP, Cambridge, 1996).

\bibitem{LeSoYu94}
G. M. Levin, M. L. Sodin, and P. Yuditskii,
{\em Ruelle operators with rational weights for Julia sets},
J. Anal. Math. {\bf 63}, 303 (1994).

\bibitem{Mart_BLMS83}
N. Martin,
{\em On finite Blaschke products whose restrictions to the unit circle are exact endomorphisms},
Bull.\ Lond.\ Math.\ Soc. {\bf 15}, 343 (1983).

\bibitem{Ni23}
A. Ni,
{\em Fast adjoint algorithm for linear responses of hyperbolic chaos},
SIAM J.\ Appl.\ Dyn.\ Syst. {\bf 22}, 2792 (2023).

\bibitem{Noet_PD19}
F. Noethen,
{\em A projector-based convergence proof of the Ginelli algorithm
for covariant Lyapunov vectors},
Physica D {\bf 396}, 18 (2019).

\bibitem{PoSe_NONL23}
M. Pollicott and B. Sewell, 
{\em Explicit examples of resonances for Anosov maps of the torus}, 
Nonlinearity {\bf 36}, 110 (2023). 

\bibitem{PuSh_ETDS08}
E. Pujals and M. Shub, {\em Dynamics of two-dimensional Blaschke products},
  Ergod.\ Th.\ \& Dyn.\ Sys. {\bf 20}, 575 (2008). 

\bibitem{Ruel_CMP97}
D. Ruelle, {\em Differentiation of SRB states},
Commun.\ Math.\ Phys. {\bf 187}, 227 (1997).

\bibitem{SlBaJu_NONL13}
J. Slipantschuk, O.F. Bandtlow, and W. Just,
{\em Analytic expanding circle maps with explicit spectra},
Nonlinearity {\bf 26},  3231  (2013).

\bibitem{SlBaJu_NONL17}
J. Slipantschuk, O.F. Bandtlow, and W. Just, {\em Complete spectral data for
  analytic Anosov maps of the torus}, Nonlinearity {\bf 30},  2667  (2017).

\bibitem{SlBaJu22}
J. Slipantschuk, O.F. Bandtlow, and W. Just,
{\em Resonances for rational Anosov maps on the torus},
(2022) Preprint: arXiv:2211.05925.

\bibitem{StHeLi_PRL85}
J. Stavans, F. Heslot, and A. Libchaber, {\em Fixed winding number and the
  quasiperiodic route to chaos in a convective fluid}, Phys.\ Rev.\ Lett. {\bf
  55},  596  (1985).

\end{thebibliography}
\end{document}